\documentclass{article}

\def\ra{\rightarrow}

\def\ss{\subseteq}

\def\Re{\hbox{\rm Re}\,}

\def\dbar{\overline{\partial}}

 \def\HollowBox #1#2{{\dimen0=#1 \advance\dimen0 by -#2       
       \dimen1=#1 \advance\dimen1 by #2                       
        \vrule height #1 depth #2 width #2                    
        \vrule height 0pt depth #2 width #1                   
        \llap{\vrule height #1 depth -\dimen0 width \dimen1}%
       \hskip -#2                                             
       \vrule height #1 depth #2 width #2}}                   
              
\font\teneufm=eufm10
\font\seveneufm=eufm7
\font\fiveeufm=eufm5
\newfam\eufmfam
\textfont\eufmfam=\teneufm
\scriptfont\eufmfam=\seveneufm
\scriptscriptfont\eufmfam=\fiveeufm

\newfam\msbfam
\font\tenmsb=msbm10  \textfont\msbfam=\tenmsb
\font\sevenmsb=msbm7  \scriptfont\msbfam=\sevenmsb
\font\fivemsb=msbm5    \scriptscriptfont\msbfam=\fivemsb
\def\Bbb{\fam\msbfam \tenmsb}

\def\RR{{\Bbb R}}
\def\CC{{\Bbb C}}

\newtheorem{theorem}{Theorem}

\makeindex

\usepackage{graphicx}

\usepackage{amsmath}

\begin{document}

\begin{center}
\huge \bf
Deformations of Strongly Pseudoconvex Domains
\end{center}
\vspace*{.12in}

\begin{center}
\large Steven G. Krantz\footnote{Author supported in part
by the National Science Foundation and by the Dean of the Graduate
School at Washington University.}\footnote{{\bf MR Classification
Numbers:}  32T15, 32T05, 32E05, 32E40.}
\end{center}
\vspace*{.15in}

\begin{center}
\today
\end{center}
\vspace*{.2in}

\begin{quotation}
{\bf Abstract:} \sl
We show that two smoothly bounded, strongly pseudoconvex domains which
are diffeomorphic may be smoothly deformed into each other, with all intermediate
domains being strongly pseudoconvex.  This result relates to Lempert's ideas
about Kobayashi extremal discs, and also has intrinsic interest.
\end{quotation}
\vspace*{.25in}

\setcounter{section}{-1}

\section{Introduction}

Ever since the solution of the Levi problem in the 1940s and 1950s, it has been
a matter of central importance to understand the geometry of pseudoconvex domains.
This investigation has many aspects, including topological features and analytic
features.

One question that has received little attention is that of deforming one pseudoconvex domain
into another.  This matter is subtle.  A form of the question comes up at the end of
Lempert's seminal paper [LEM], where he deforms strongly {\it convex} domains.  Of course
that is a relatively easy matter, but it begs the question (if one wants to generalize Lempert's
results to strongly {\it pseudoconvex} domains) of deforming more general types of domains.
See [KRA2] for an investigation of this type of generalization.

In the present paper we explore such deformation ideas, and we prove a positive
deformation result for smoothly bounded, strongly pseudoconvex domains.  One interesting
aspect of our approach is that we make decisive use of the Forn\ae ss embedding theorem [FOR].
This may be the first actual application of Forn\ae ss's theorem.  Of course Forn\ae ss's theorem
applies only to strongly pseudoconvex domains, and there is no hope of adapting the techniques
presented here to a more general class of domains.  It is plausible that various natural classes
of Levi geometry should be preserved under smooth deformation, but the techniques for proving
a very general result do not seem to be available at this time.

It is a pleasure to thank Jiye Yu for helpful conversations on the subject matter of this
paper.

\section{Principal Results}

Our main theorem is this:

\begin{theorem} \sl 
Let $\Omega_{-1}$, $\Omega_1$ be smoothly
bounded, strongly pseudoconvex domains in $\CC^n$. Write
$\Omega_j = \{z \in \CC^n: \rho_j(z) < 0\}$, $j = -1, 1$, where
$\rho_j$ is a defining function for $\Omega_j$ (see [KRA1] for
this concept).	 Assume that $\overline{\Omega}_{-1}$ and $\overline{\Omega}_1$
are diffeomorphic.

Then there is a smooth deformation $\rho(z,t)$ on $\CC^n \times [-1,1]$ so that,
for each $t \in [-1,1]$, the
function $\rho_t(z) \equiv \rho(z,t)$ is a smooth defining function.  
Also $\rho(z,0) = \rho_{-1}(z)$ and $\rho(z,1) = \rho_1(z)$.   Moreover, for each $t \in [-1,1]$,
the domain $\Omega_t \equiv \{z \in \CC^n: \rho_t(z) < 0\}$ is smoothly bounded and strongly
pseudoconvex.	Finally, there is a $c > 0$ so that all the eigenvalues of the Levi form at
any boundary point of any $\Omega_t$ are not less than $c$.
\end{theorem}

It should be stressed that we do {\it not} conclude that the domains $\Omega_{-1}$ and $\Omega_1$ 
are biholomorphic.  This is impossible in the strongest sense---see [GRK] and [BSW].
We are only setting up a smooth deformation of one domain into the other.  What is interesting
about the proof is that we simultaneously Forn\ae ss-embed the two domains $\Omega_{-1}$ and $\Omega_1$
into a single strongly convex domain $W$, perform the deformation in $W$, and then pull it back.
Along the way, an uncountably infinite family of $\overline{\partial}$ problems must be solved.

\section{Proof of the Theorem}

According to the Forn\ae ss embedding theorem [FOR], there is, for $j
= -1, 1$, a strongly convex domain $W_j \ss \CC^{N_j}$ (with,
in general $N_j > \, > n$), a neighborhood $U_j$ of $\overline{\Omega}_j$, 
and a univalent holomorphic mapping
$\Phi_j : U_j \ra \CC^{N_j}$ with
$\Phi_j(\Omega_j) \ss W_j$, $\Phi_j(\overline{\Omega}_j) \ss
\overline{W}_j$, $\Phi_j(U_j \setminus \overline{\Omega}_j) \subseteq \CC^{N_j} \setminus \overline{W}_j$,
and so that $\Phi_j(\overline{\Omega}_j)$ is
transversal to $\partial W_j$ where they meet.

We may assume that $N_1 = N_{-1}$.  Let $N = 2N_1 + 1$.  Let $\varphi: \RR \ra \RR$ be an even, strictly concave
function so that
\begin{enumerate}
\item[{\bf (1)}]  $\varphi(-1) = 0$;
\item[{\bf (2)}]  $\varphi(1) = 0$;
\item[{\bf (3)}]  $\varphi(0) = 1$;
\item[{\bf (4)}]  $\varphi$ is strictly decreasing to $-\infty$ on $(0,\infty)$;
\item[{\bf (4)}]  $- \varphi''(x) \geq K$ for some large positive $K > 0$ and all $x \in \RR$.
\end{enumerate}
Also let $\lambda: \RR \ra \RR$ be a smooth function satisfying
\begin{enumerate}
\item $\lambda$ is smooth;
\item $0 \leq \lambda(x) \leq 1$ for all $x$;
\item $\lambda(x) \equiv 0$ when $x \leq -1$;
\item $\lambda(x) \equiv 1$ when $x \geq 1$;
\item $\lambda$ is monotone increasing for $-1 < x < 1$;
\end{enumerate}

Let $\mu_{-1}$ be a defining function for $W_{-1}$ and $\mu_1$ a defining function for $W_1$.

If $z \in \CC^N$, then let us write
$$
z = (z_1, z_2, \dots, z_N) = (z'; z^*) \, ,
$$
where
$$
z' = (z_1, z_2, \dots, z_{N_1}) 
$$
and
$$
z^* = (z_{N_1 + 1}, z_{N_1 + 2}, \dots, z_N) \, .
$$

Now we consider in $\CC^N$ the domain with defining function
$$
\rho(z) = (1 - \lambda(\Re z_N)) \mu_{-1}(z') + \lambda(\Re z_N) \mu_1(z') - \varphi(|z^*|^2) \, .
$$
We set
$$
W = \{z \in \CC^N : \rho(z) < 0\} \, .
$$
We claim that, in a natural sense, $\Omega_{-1}$ and $\Omega_1$ are Forn\ae ss-embedded into $W$.
Let us see why.

First, $\rho$ is strictly convex.  The pure second derivatives of $\rho$ in $z'$ are clearly under control
and positive.  The pure second derivatives in $z^*$ are controlled by $- \varphi''$ provided that $K$ is large
enough.  Also the mixed second derivatives are controlled by $-\varphi''$ provided
that $K$ is large enough.

Second, $\rho(z'; (0, 0, \dots, -1 + 0i)) = \mu_{-1}(z')$ so that $\{z: z^* = (0,0, \dots, -1+0i), \rho(z) < 0\}$ is simply
a copy of $W_{-1}$.  And $\rho(z'; (0, 0, \dots, 1 + 0i)) = \mu_1(z')$, so that $\{z: z^* = (0,0,\dots,1+0i), \rho(z) < 0\}$ is
a copy of $W_1$.  Thus 
$$
\widetilde{\Phi}_{-1}(w) \equiv (\Phi_{-1}(w); (0, 0, \dots, -1+0i))
$$
embeds $\Omega_{-1}$ into $W$.  And
$$
\widetilde{\Phi}_1(w) \equiv (\Phi_1(w); (0,0,\dots, 1+0i))
$$
embeds $\Omega_1$ into $W$.   As a result, we have both $\Omega_{-1}$ and $\Omega_1$ embedded,
in the fashion of the Forn\ae ss embedding theorem, into the single strongly convex
domain $W$.

Assume without loss of generality that 0 lies in $W$.  Now we may define, for $z \in \CC^n$, and $-1 \leq t \leq 1$, 
a smooth function $\omega(z,t)$ so that $\omega(z,-1)$ is a defining
function for $W \cap \{z: z^* = (0,0,\dots,-1+0i)\}$, $\omega(z, 1)$ is a defining
function for $W \cap \{z: z^* = (0,0,\dots,1+0i)\}$, and, in general, $\omega(z,t)$ is
a defining function for $W \cap \{z: z^* = (0,0,\dots, t+ i0)\}$, $-1 \leq t \leq 1$.  Let $S_t = \{z \in W: z^* = (0,0,\dots,t+0i)\}$.
Then, for $t$ near $-1$, we may map $S_{-1}$ to $S_t$ by sending a point $z$ in $S_{-1}$ that is distance
$\delta$ from the boundary to a point $\widetilde{z}_t$ in $S_t$ that is distance $\delta$ from the boundary
and so that the line through the origin and $\widetilde{z}_t$ has the same projection to the set
$\{z^* = 0\}$ as the line through the origin and $z$.   This is a smooth mapping for $z$ near
the boundary, and we may easily interpolate it to a smooth mapping on all of $S_{-1}$.  Call the mapping $\pi_t$.
Then, if $\widetilde{\Phi}_{-1}$ is the Forn\ae ss embedding of $\Omega_{-1}$ into $S_{-1}$ and $\widetilde{\Phi}_{-1}^{-1}$ its inverse
defined on the image of $\widetilde{\Phi}_{-1}$, then we have
a pseudo-inverse-embedding $\widetilde{\Phi}_{-1}^{-1} \circ \pi_t^{-1}$ of $S_t$ into $\CC^n$.  We call this a pseudo-inverse-embedding
because it will not be holomorphic.  But we can solve the $\overline{\partial}$ equation (see [HEN] and [SIU])
$$
\overline{\partial} h = - \overline{\partial} \biggl [ \widetilde{\Phi}_{-1}^{-1} \circ\pi_t^{-1} \biggr ]
$$
to find a function $h$ (which will have small $C^2$ norm because $\overline{\partial} \pi_t^{-1}$ has
small $C^2$ norm) so that $\tau_t \equiv \widetilde{\Phi}_{-1}^{-1} \circ \pi_t^{-1} + h$ is a holomorphic embedding of $S_t$ into $\CC^n$.
We think of the image of $\tau_t$ as a perturbation $\Omega_t$ of $\Omega_{-1}$.  We can continue
to incrementally increase the value of $t$ and create additional deformations of $\Omega_{-1}$ as $t$ increases.
Note here that it is propitious to use the Henkin solution of the $\overline{\partial}$ problem---see [HEN].  For it
is known [GRK] to provide smoothly varying solutions when the data is varied smoothly.

At the same time, we could begin at $\Omega_1$ and deform in the same fashion, decreasing values of $t$.  This
incrementally creates deformations of $\Omega_1$.  When the values of $t$ from above and the values of $t$ from
below meet (say at $t=0$), they must of course have data sets $S_t$ that are close together in the smooth topology.  So the
deformations constructed in $\CC^n$ will also be close together.   In sum, the two streams of deformed domains
give rise to a deformation of $\Omega_{-1}$ to $\Omega_1$.  Finally, it is clear by inspection that the eigenvalues
of the Levi form for $\Omega_t$ are pullbacks of eigenvalues of the Levi form for $W$.  So the eigenvalues
of the Levi form for all the $\Omega_t$ are uniformly bounded from 0.

\section{Concluding Remarks}

It would be a matter of some interest to solve the problem treated here in
the algebraic category: If $\Omega_{-1}$ and $\Omega_1$ are bounded
domains in $\CC^n$ with strongly pseudoconvex, algebraic boundaries, then
can one be deformed into the other with all intermediate domains being
algebraic and strongly pseudoconvex? Likewise, one would like to solve
this problem in the real analytic category: If $\Omega_{-1}$ and
$\Omega_1$ are bounded domains in $\CC^n$ with strongly pseudoconvex, real
analytic boundaries, then can one be deformed into the other with all
intermediate domains being real analytic and strongly pseudoconvex?
Unfortunately the methods of the present paper cannot be applied to either
of these problems. Our constructions are strictly real-variable in nature.
The methods used in Grauert's solution of the Levi problem, and in his
embedding theorem for real analytic manifolds (see [GRA]), may be useful in studying
some of these new questions.

It also would be worthwhile to study these deformation questions for other types of domains, such
as pseudoconvex, finite type domains.  Certainly the Forn\ae ss embedding theorem is not
true for such domains, and we do not know how to attack such a problem at this time.

We hope to explore these new questions in a future paper.  										    

\newpage

\noindent {\Large \sc References}
\bigskip  \\

\begin{enumerate}

\item[{\bf [BSW]}] D. Burns, S. Shnider, and R. O. Wells, On
deformations of strictly pseudoconvex domains, {\em Invent.
Math.} 46(1978), 237--253.

\item[{\bf [FOR]}] J. E. Forn\ae ss, Strictly pseudoconvex
domains in convex domains, {\em Am. J. Math.} 98(1976),
529--569.

\item[{\bf [GRA]}]  H. Grauert, On Levi's problem and the imbedding of real-analytic manifolds,  
{\it Ann.\ of Math.} 68(1958), 460--472.

\item[{\bf [GRK]}] R. E. Greene and S. G. Krantz, Deformation
of complex structures, estimates for the $\dbar$ equation, and
stability of the Bergman kernel, {\em Adv. Math.} 43(1982),
1--86.

\item[{\bf [HEN]}] G. M. Henkin, Integral representations of
functions holomorphic in strictly pseudoconvex domains and
applications to the $\dbar$ problem, {\em Mat. Sb.} 82(124),
300-308 (1970); {\em Math. U.S.S.R. Sb.} 11(1970), 273-281.
																
\item[{\bf [KRA1]}] S. G. Krantz, {\it Function Theory of
Several Complex Variables}, $2^{\rm nd}$ ed., American
Mathematical Society, Providence, RI, 2001.

\item[{\bf [KRA2]}] S. G. Krantz, The Kobayashi metric,
extremal discs, and biholomorphic mappings, {\it Complex
Variables and Elliptic Equations}, to appear.

\item[{\bf [LEM]}] L. Lempert, La metrique Kobayashi et
la representation des domains sur la boule, {\em Bull. Soc.
Math. France} 109(1981), 427--474.

\item[{\bf [SIU]}] Y.-T. Siu, The $\overline{\partial}$ problem
with uniform bounds on derivatives, {\it Math.\ Ann.}
207(1974), 163--176.

\end{enumerate}
\vspace*{.95in}

\noindent \begin{quote}
Department of Mathematics \\
Washington University in St. Louis \\
St.\ Louis, Missouri 63130 \\ 
{\tt sk@math.wustl.edu}
\end{quote}

\end{document}